\documentclass{amsart}
\usepackage{amssymb}
\usepackage{latexsym}
\usepackage{amsmath}
\usepackage{euscript}
\usepackage{graphics}
\usepackage[all]{xy}

\newcommand{\be}{\begin{equation}}
\newcommand{\ee}{\end{equation}}
\newcommand{\ba}{\begin{eqnarray}}
\newcommand{\ea}{\end{eqnarray}}
\newcommand{\baa}{\begin{eqnarray*}}
\newcommand{\eaa}{\end{eqnarray*}}
\newcommand{\bb}{\begin{thebibliography}}
\newcommand{\eb}{\end{thebibliography}}

\newcommand{\bi}[1]{\bibitem{#1}}
\newcommand{\lab}[1]{\label{#1}}
\newcommand{\re}[1]{(\ref{#1})}

\newcounter{my}
\newcommand{\he}%
   {\stepcounter{equation}\setcounter{my}%
   {\value{equation}}\setcounter{equation}0%
   }%
\newcommand{\she}%
   {\setcounter{equation}{\value{my}}%
    }%

\renewcommand\t{\tilde}

\theoremstyle{definition}

\numberwithin{equation}{section}

\begin{document}

\title[A limit $q=-1$ for the big q-Jacobi polynomials]
{A limit $q=-1$ for the big q-Jacobi polynomials}

\author{Luc Vinet}
\author{Alexei Zhedanov}

\address{Centre de recherches math\'ematiques
Universite de Montr\'eal, P.O. Box 6128, Centre-ville Station,
Montr\'eal (Qu\'ebec), H3C 3J7}

\address{Institute for Physics and Engineering\\
R.Luxemburg str. 72 \\
83114 Donetsk, Ukraine \\}

%\email{zhedanov@yahoo.com}

\date{\today}

\begin{abstract}
We study a new family of "classical" orthogonal polynomials, here
called big -1 Jacobi polynomials, which satisfy (apart from a
3-term recurrence relation) an eigenvalue problem  with
differential operators of Dunkl-type. These polynomials can be
obtained from the big $q$-Jacobi polynomials in the limit $q \to
-1$. An explicit expression of these polynomials in terms of
Gauss' hypergeometric functions is found. The big -1 Jacobi
polynomials are orthogonal on the union of two symmetric intervals
of the real axis. We show that the big -1 Jacobi polynomials can
be obtained from the Bannai-Ito polynomials when the orthogonality
support is extended to an infinite number of points. We further
indicate that these polynomials provide a nontrivial realization
of the Askey-Wilson algebra for $q \to -1$.
\end{abstract}


\maketitle

\section{Introduction}
\setcounter{equation}{0}
 We constructed in \cite{VZ_little} a system of "classical`` orthogonal polynomials $P_n(x)$
 containing two real parameters $\alpha,\beta$ and corresponding to the limit $q \to -1$ of the little q-Jacobi polynomials.
 By ''classical'' we mean that these polynomials satisfy (apart from a 3-term recurrence relation) a nontrivial eigenvalue equation of the form
\be
L P_n(x) = \lambda_n P_n(x). \lab{eig_P} \ee
The novelty lies in the fact that $L$ is a differential-difference operator of special type.
Namely, $L$ is a linear operator which is of first order in the derivative operator $\partial_x$
and contains also the reflection operator $R$ which acts as $Rf(x)=f(-x)$. Roughly speaking,
one can say that $L$ belongs to the class of Dunkl operators \cite{Dunkl} which contain both the operators $\partial_x$ and $R$.
Nevertheless, the operator $L$ differs from the standard Dunkl operators in a fundamental way.
Indeed, $L$ preserves the linear space of polynomials of any given maximal degree.
This basic property allows to construct a complete system of polynomials $P_n(x), \: n=0,1,2,\dots$ as eigenfunctions of the operator $L$.

Guided by the $q \to -1$ limit of the little $q$-Jacobi
polynomials, we derived in \cite{VZ_little} an explicit expression
of the polynomials $P_n(x)$ in terms of Gauss' hypergeometric
functions. We also found explicitly the recurrence coefficients
and showed that the polynomials $P_n(x)$ are orthogonal on the
interval $[-1,1]$ with a weight function related to the weight
function of the generalized Jacobi polynomials \cite{Chi_Chi}. We
also proved that they admit the Dunkl classical property
\cite{Cheikh} and further demonstrated that the operator $L$
together with the multiplication operator $x$ form a special case
of the Askey-Wilson algebra $AW(3)$ \cite{Zhe} corresponding to
the parameter $q=-1$.

In this paper we construct similarly, a new family of "classical`` orthogonal polynomials which are
obtained as a nontrivial limit of the  big q-Jacobi polynomials when $q \to -1$. We will call them ''big -1 Jacobi polynomials``

In contrast to the little -1 Jacobi polynomials,  the big -1
Jacobi polynomials contain 3 real parameters $\alpha,\beta,c$.
This leads to more complicated formulas for the recurrence
coefficients as well as for the explicit expression in terms of
the Gauss hypergeometric function. Moreover, in contrast to the
little -1 Jacobi polynomials the big -1 Jacobi polynomials are
orthogonal on the union of the {\it two intervals} $[-1,-c]$ and
$[c,1]$ (it is assumed that $0<c<1$). When $c=0$ these intervals
connect into one interval $[-1,1]$. This corresponds to the
degeneration of the big -1 Jacobi polynomials into the little -1
Jacobi polynomials

The fundamental ''classical`` property \re{eig_P} holds for the big -1 Jacobi polynomials as well.
The operator $L$ is again a first order differential operator of Dunkl type which preserves the space of polynomials.
This means that both little and big -1 Jacobi polynomials provide two ''missing'' families of classical orthogonal
polynomials which should be included into the Askey table as special cases.

Limit cases $q \to -1$ for q-ultraspherical polynomials were
considered by Askey and Ismail \cite{AI2}, \cite{AI}. Bannai and
Ito found a nontrivial example of explicit polynomials
corresponding to the limit $q \to -1$ in the Askey scheme was
proposed by \cite{BI}, \cite{Ter2}. The Bannai-Ito polynomials are
orthogonal on a finite set of $N+1$ points. We show that the big
-1 Jacobi polynomials are obtained in the limit $N \to \infty$
from the Bannai-Ito polynomials.

We also demonstrate that the big -1 Jacobi polynomials provide a
convenient realization of the $AW(3)$ algebra for $q=-1$.

\section{Big q-Jacobi polynomials in the limit $q=-1$}
\setcounter{equation}{0} The big q-Jacobi polynomials
$P_n(x;a,b,c)$ were introduced by W. Hahn \cite{Hahn} in 1949.
Andrews and Askey found explicit orthogonal relation for these
polynomials in \cite{AA}. Implicitly, the big q-Jacobi polynomials
are also contained in the Bannai-Ito scheme of dual systems of
orthogonal polynomials as an infinite dimension analogue of the
q-Racah polynomials \cite{BI}. These polynomials depend on 3
parameters $a,b,c$ and are defined by the following 3-term
recurrence relation (for brevity, we will sometimes omit the
dependence on the parameters $a,b,c$): \be P_{n+1}(x) + b_n P_n(x)
+ u_n P_{n-1}(x) = x P_n(x), \lab{3-term} \ee where
$$
u_n = A_{n-1}C_n, \quad b_n = 1-A_n-C_n
$$
with \be A_n =
\frac{(1-aq^{n+1})(1-abq^{n+1})(1-cq^{n+1})}{(1-abq^{2n+1})(1-abq^{2n+2})},
\quad C_n = -acq^{n+1} \:
\frac{(1-q^{n})(1-abc^{-1}q^{n})(1-bq^{n})}{(1-abq^{2n+1})(1-abq^{2n})}
\lab{AC_qJ} \ee In terms of basic hypergeometric functions
\cite{KS}, \cite{KLS} they are given by \be P_n(x;a,b,c)= \kappa_n
{_3}\varphi_2 \left( {q^{-n}, abq^{n+1},x  \atop  aq, cq} \Big  |
q; q   \right) \lab{P_hyp} \ee where the coefficient $\kappa_n$
ensures that $P_n(x)$ is monic: $P_n(x) = x^n + O(x^{n-1})$. We
shall not need the explicit expression of $\kappa_n$ in the
following.

The big q-Jacobi polynomials satisfy the eigenvalue equation \cite{KS}, \cite{KLS}
\be
L P_n(x) = \lambda_n P_n(x), \quad \lambda_n = (q^{-n}-1)(1-abq^{n+1})  \lab{EVP_qJ} \ee
where the operator $L$ is
\be
L f(x)= B(x) (f(xq)-f(x)) + D(x)(f(xq^{-1}) - f(x)) \lab{L_op} \ee
with
\be
B(x)=\frac{aq(x-1)(bx-c)}{x^2}, \quad D(x) = \frac{(x-aq)(x-cq)}{x^2} \lab{BD_qJ} \ee
The orthogonality relation is
\be
\int_{cq}^{aq} w(x) P_n(x) P_m(x) d_q x = h_n \delta_{nm}, \quad h_n = u_1 u_2 \dots u_n, \lab{ort_qJ} \ee
with the q-integral defined as \cite{KS}, \cite{KLS}
$$
\int_{cq}^{aq} f(x) d_q x= aq(1-q)\sum_{s=0}^{\infty} f(aq^{s+1})q^s -  cq(1-q)\sum_{s=0}^{\infty} f(cq^{s+1})q^s
$$
and the weight function
\be
w(x) =g \frac{(a^{-1}x;q)_{\infty} (c^{-1}x;q)_{\infty}}{(x;q)_{\infty} (bc^{-1}x;q)_{\infty}}, \lab{w_big} \ee
where
$$
(a;q)_s=(1-a)(1-aq) \dots(1-aq^{s-1})
$$
is the shifted q-factorial \cite{KLS} and $(a;q)_{\infty} = \lim_{s \to \infty} (a;q)_s$
(In \re{w_big}, $g$ is a normalization factor which is not essential for our considerations).

Consider the operator $(q+1)^{-1}L$, where the operator $L$ is
defined by \re{L_op}. Put \be q=-\exp(\epsilon),\:
a=-\exp(\epsilon \alpha), \: b=-\exp(\epsilon \beta) \lab{qab_lim}
\ee and take the limit $\epsilon \to 0$ which corresponds to the
limit $q \to -1$. It is not difficult to verify that the limit
does exist and that we have \be L_0 = \lim_{q \to -1} (q+1)^{-1}L
= g_0(x)(R-I) + g_1(x) \partial_x R, \lab{L_0} \ee where \be
g_0(x)= \frac{(\alpha+\beta+1)x^2 +(c\alpha -\beta)x + c}{x^{2}},
\quad g_1(x)=\frac{2(x-1)(x+c)}{x}. \lab{g_01} \ee The operator
$I$ is the identity operator and $R$ is the reflection operator
$Rf(x) = f(-x)$.

Equivalently, the operator $L_0$ can be presented through its action on $f(x)$:
\be
L_0f(x) = g_0(x)(f(-x)-f(x)) - g_1(x)f'(-x) \lab{L_0f} \ee
On monomials $x^n$ the operator $L_0$ acts as follows.

For $n$ even,
\be
L_0 x^n = 4n(x-1)(x+c) x^{n-2}. \lab{L_x_even} \ee

For $n$ odd
\be
L_0 x^n = -2(\alpha+\beta+n+1)x^n + 2(\beta-c\alpha+n-c)x^{n-1} + 2(n-1)cx^{n-2} \lab{L_x_odd} \ee
In any case, the operator $L_0$ is lower triangular, with 3 diagonals, in the basis $x^n$:
\be
L x^n = \xi_n x^n +\eta_n x^{n-1} + \zeta_n x^{n-2} \lab{L-3-diag_x} \ee with the coefficients $\xi_n,\eta_n,\zeta_n$
straightforwardly obtained from \re{L_x_even}, \re{L_x_odd}.
It is easily seen that the operator $L_0$ preserves the linear space of polynomials of any fixed dimension. Hence for every $n=0,1,2,\dots$
there are monic polynomials eigenfunctions $P_n^{(-1)}(x)=x^n + O(x^{n-1})$ of the operator $L_0$

This eigenvalue equation is obtained as the  $q \to -1$ limit of
the eigenvalue equation \re{EVP_qJ}: \be L_0 P_n^{(-1)}(x) =
\lambda_n \: P_n^{(-1)}(x),  \lab{eigen_P-1} \ee where \be
\lambda_n = \left\{ {2n, \quad n \quad \mbox{even}    \atop
-2(\alpha+\beta+n+1), \quad n \quad \mbox{odd}} \right .
\lab{lam-1} \ee Consider the limit $q \to -1$ for the recurrence
coefficients. Assuming \re{qab_lim}, we have \be A_n^{(-1)} =
\lim_{\epsilon \to 0} A_n = \left\{  {\frac{(c+1)
(\alpha+n+1)}{\alpha+\beta+2n+2}, \quad n \quad \mbox{even} \atop
\frac{(1-c) (\alpha+\beta+n+1)}{\alpha+\beta+2n+2}, \quad n \quad
\mbox{odd} } \right . \lab{A_lim} \ee and \be C_n^{(-1)} =
\lim_{\epsilon \to 0} C_n = \left\{
{\frac{(1-c)n}{\alpha+\beta+2n}, \quad n \quad \mbox{even} \atop
\frac{(1+c) (\beta+n)}{\alpha+\beta+2n}, \quad n \quad \mbox{odd}
} \right . \lab{C_lim} \ee Hence for the recurrence coefficients
we have \be u_n^{(-1)} = \lim_{\epsilon \to 0} {A_{n-1}C_n} =
\left\{  {\frac{(1-c)^2 n(\alpha+\beta+n)}{(\alpha+\beta+2n)^2},
\quad n \quad \mbox{even} \atop  \frac{(1+c)^2
(\alpha+n)(\beta+n)}{(\alpha+\beta+2n)^2}, \quad n \quad
\mbox{odd} } \right . \lab{u_lim} \ee and \be b_n^{(-1)} =
\lim_{\epsilon \to 0} {1-A_{n}-C_n} = \left\{  {-c
+\frac{(c-1)n}{\alpha+\beta+2n}
+\frac{(1+c)(\beta+n+1)}{\alpha+\beta+2n+2}, \quad n \quad
\mbox{even} \atop  c + \frac{(1-c)(n+1)}{\alpha+\beta+2n+2} -
\frac{(c+1)(\beta+n)}{\alpha+\beta+2n}, \quad n \quad \mbox{odd} }
\right . \lab{b_lim} \ee The polynomials $P_n^{(-1)}(x)$ satisfy
the 3-term recurrence relation \be P_{n+1}^{(-1)}(x) + b_n^{(-1)}
P_n^{(-1)}(x) + u_n^{(-1)} P_{n-1}^{(-1)}(x) = xP_n^{(-1)}(x).
\lab{rec-1} \ee For any real $c \ne 1$ and real $\alpha, \beta$
satisfying the restriction $\alpha>-1, \: \beta>-1$, the
recurrence coefficients $b_n^{(-1)}$ are real and the recurrence
coefficients $u_n$ are positive. This means that the polynomials
$P_n^{(-1)}(x)$ are positive definite orthogonal polynomials.

Let us consider expression \re{P_hyp} in details,
\be
P_n(x)= \kappa_n \: \sum_{s=0}^n \frac{(q^{-n};q)_s    (abq^{n+1};q)_s (x;q)_s}{(q;q)_s (aq;q)_s (cq;q)_s} q^s, \lab{P_q_hyp} \ee
In the limit $q \to -1$ it is easy to obtain that
$$
\frac{(x;q)_s}{(cq;q)_s} = \left\{   \left(\frac{1-x^2}{1-c^2}\right)^{s/2}, \quad s \quad \mbox{even} \atop  \frac{1-x}{1+c}  \left(\frac{1-x^2}{1-c^2}\right)^{(s-1)/2}, \quad s \quad \mbox{odd}    \right .
$$
Hence, in the limit $q \to -1$ the sum \re{P_q_hyp} is divided in two parts.
The first part is an even polynomial with respect to $x$, i.e. $p(x^2)$, where $p(x)$ is a polynomial.
The second part will have the form $(1-x) q(x^2)$  with another polynomial $q(x)$. Simple calculations lead to the following formulas.

If $n$ is even
\be
P_n^{(-1)}(x) = \kappa_n \left[ {_2}F_1\left( {-\frac{n}{2}, \frac{n+\alpha+\beta+2}{2} \atop \frac{\alpha+1}{2}} \left |    \frac{1-x^2}{1-c^2}   \right)   \right .     +   \frac{n(1-x)}{(1+c)(\alpha+1)}  {_2}F_1\left( {1-\frac{n}{2}, \frac{n+\alpha+\beta+2}{2} \atop \frac{\alpha+3}{2}} \left |    \frac{1-x^2}{1-c^2}   \right)   \right . \right]  \lab{even_hyp} \ee

If $n$ is odd \be P_n^{(-1)}(x) = \kappa_n \left[ {_2}F_1\left(
{-\frac{n-1}{2}, \frac{n+\alpha+\beta+1}{2} \atop
\frac{\alpha+1}{2}} \left |    \frac{1-x^2}{1-c^2}   \right)
\right .     -   \frac{(\alpha+\beta+n+1)(1-x)}{(1+c)(\alpha+1)}
{_2}F_1\left( {-\frac{n-1}{2}, \frac{n+\alpha+\beta+3}{2} \atop
\frac{\alpha+3}{2}} \left |    \frac{1-x^2}{1-c^2}   \right)
\right . \right]  \lab{odd_hyp} \ee The normalization coefficient
is given by \be \kappa_n = \left\{   { \frac{(1-c^2)^{n/2}
((\alpha+1)/2)_{n/2}  }{((n+\alpha+\beta+2)/2)_{n/2}}, \quad n
\quad \mbox{even} \atop (1+c)\frac{(1-c^2)^{(n-1)/2}
((\alpha+1)/2)_{(n+1)/2}  }{((n+\alpha+\beta+1)/2)_{(n+1)/2}},
\quad n \quad \mbox{odd}}   \right . \lab{kappa_n} \ee The
remaining problem is to find the orthogonality relation and the
corresponding weight function $w(x)$ for the big -1 Jacobi
polynomials. Of course, this could be done directly from the known
orthogonality relation for the big q-Jacobi polynomials by taking
the limit $q \to -1$. However it is more instructive to derive the
weight function using the method of polynomial mappings
\cite{GVA}, \cite{MP}. This method will allow to find nontrivial
relations between the big -1 Jacobi polynomials and  the ordinary
Jacobi polynomials. This will explain the origin of the rather
``strange`` expressions \re{even_hyp} and \re{odd_hyp}.

\section{Polynomial systems and the Christoffel transform}
\setcounter{equation}{0}
In this section we consider a scheme allowing to obtain a new family of orthogonal polynomial starting from two sets of orthogonal
polynomials related by the  Christoffel transform. This scheme is a simple generalization of the well known Chihara
method for constructing symmetric orthogonal polynomials from a pair of orthogonal polynomials and their kernel partner \cite{Chi}.
It is also very close to the scheme proposed by Marcell\'an and Petronilho in \cite{MP}.

Let $P_n(x), \; n=0,1,2, \dots $ be a set of monic orthogonal polynomials satisfying the recurrence relation
\be
P_{n+1}(x) + b_n P_n(x) + u_n P_{n-1}(x) = x P_n(x). \lab{rec_P} \ee
Consider a partner family of orthogonal polynomials $Q_n(x)$ related to $P_n(x)$ by the Christoffel transform \cite{Sz}
\be
Q_n(x) = \frac{P_{n+1}(x) - A_n P_n(x)}{x- \nu^2 }, \lab{Q_P_CT} \ee
where $\nu$ is a real parameter and  $A_n= P_{n+1}(\nu^2)/P_n(\nu^2)$.

If the polynomials $P_n(x)$ are monic orthogonal with respect to the linear functional $\sigma$
$$
\langle \sigma, P_n(x) P_m(x) \rangle =0, \quad n \ne m
$$
then the polynomials $Q_n(x)$ are monic orthogonal with respect to the functional $\t \sigma = (x-\nu^2) \sigma$, i.e. \cite{Sz}
$$
\langle \sigma, (x-\nu^2)Q_n(x) Q_m(x) \rangle =0, \quad n \ne m
$$
The polynomials $P_n(x)$ are expressed in terms of the polynomials $Q_n(x)$ via the Geronimus transform \cite{ZhS}
\be
P_n(x) = Q_n(x) - B_n Q_{n-1}(x), \lab{PQ_GT} \ee
where the coefficients $B_n$ are related to $A_n$ and the recurrence coefficients by the formulas
\be
u_n = B_n A_{n-1}, \quad b_n = -A_n-B_n+\nu^2. \lab{ub_AB} \ee
Now, starting from a pair of polynomials $P_n(x), Q_n(x)$ we can construct another family of orthogonal polynomial  $R_n(x)$ by proceeding as follows.

For even numbers $n$, let the polynomials $R_n(x)$ be defined according to
\be
R_{2n}(x) = P_n(x^2) \lab{RP_even} \ee
and for odd numbers $n$, let
\be
R_{2n+1}(x) = (x-\nu)Q_n(x^2) \lab{RP_odd} \ee
It is obvious that for all $n=0,1,2,\dots$ the polynomials $R_n(x)$ are monic polynomials in $x$ of degree $n$.

What is more important is that the polynomials $R_n(x)$ are {\it orthogonal}, since they satisfy the 3-term recurrence relation
\be
R_{n+1}(x) + (-1)^n \nu R_n(x) + v_n R_{n-1}(x) = x R_n(x), \lab{3term_R} \ee
where
\be
v_{2n}=-B_n, \; v_{2n+1}=-A_n \lab{xi_eta_AB} \ee
This construction can also be carried out in the reverse.

Assume that the polynomials $R_n(x)$ satisfy the recurrence relation \re{3term_R} with some real parameter $\nu$
and positive coefficients $v_n$, it can easily be shown by induction that
$$
R_{2n}(x) = P_n(x^2), \quad R_{2n+1}(x) = (x-\nu) Q_n(x^2),
$$
where $P_n(x), Q_n(x)$ are monic polynomials of degree $n$.

The polynomials $R_n(x)$ are orthogonal with respect to a positive definite linear functional $\rho$:
\be
\langle \rho, R_n(x) R_m(x) \rangle = 0, \quad n \ne m \lab{ort_R} \ee
Let
$$
r_n = \langle \rho, x^n \rangle
$$
be the corresponding moments. We use the standard normalization condition $r_0=1$. It can then be proven, again by induction, that
\be
r_{2n+1} = \nu r_{2n}, \quad n=0,1,2,\dots \lab{r_mom_cond} \ee
and that the even moment $r_{2n}$ is an even monic polynomial of degree $2n$ in the argument $\nu$, i.e.
$$
r_{2n} = \nu^{2n} + n v_1 \nu^{2n-2} + \frac{n(n-1)}{2} v_1(v_1+v_2) \nu^{2n-4} + O(\nu^{2n-6}).
$$
It is directly verified that the polynomials $P_n(x)$ and $Q_n(x)$ are orthogonal as they satisfy the recurrence relations
$$
P_{n+1}(x) +(v_{2n} + v_{2n+1}+ \nu^2)P_n(x)+ v_{2n} v_{2n-1} P_{n-1}(x)=xP_n(x)
$$
and
$$
Q_{n+1}(x) +(v_{2n+2} + v_{2n+1}+ \nu^2)Q_n(x)+ v_{2n} v_{2n+1} Q_{n-1}(x)=xQ_n(x)
$$
Moreover, the polynomials $Q_n(x)$ are Christoffel transforms of the polynomials $P_n(x)$:
$$
Q_n(x) = \frac{P_{n+1}(x) + v_{2n+1}P_n(x)}{x-\nu^2}.
$$
while the polynomials $P_n(x)$ are Geronimus transforms of $Q_n(x)$:
$$
P_n(x) = Q_n(x) + v_{2n} Q_{n-1}(x)
$$
Assume that the polynomials $P_n(x)$ have moments $c_n$. Then one has a simple relation between the moments
\be
r_{2n} = c_{n}, \; r_{2n+1}= \nu c_n, \quad n=0,1,2,\dots  \lab{r_c_rel} \ee
The moments $\t c_n$ corresponding to the polynomials $Q_n(x)$ are given by
\be
\t c_n = \frac{c_{n+1}- \nu^2 c_n}{c_1- \nu^2} \lab{tc_exp} \ee
Expression \re{tc_exp} follows easily from the definition of the Christoffel transform \cite{ZhS}.

Note that in the special case $\nu=0$ we recover the well known scheme relating symmetric and non-symmetric
polynomials that has been described in details by Chihara \cite{Chi}.
In this case the polynomials $R_n(x)$ are symmetric $R_n(-x) =(-1)^n R_n(x)$
and their odd moments are zero $r_{2n+1}=0$. All the above formulas remain valid if one puts $\nu=0$.
We have thus provided a generalization of the Chihara scheme with an additional parameter $\nu$.
Note that the resulting polynomials $R_n(x)$ are no longer symmetric; they satisfy however
the simple recurrence relation \re{3term_R} and have properties very close to those of symmetric orthogonal polynomials.

In \cite{MP} a more general problem was studied with the orthogonal polynomials $R_n(x)$ defined as $R_{2n}(x)=P_n(\phi(x))$, where $\phi(x)$
is a  polynomial of second degree and $P_n(x)$ is a given system of orthogonal polynomials. Our approach corresponds to  the special case $\phi(x)=x^2$. Note that the general case of polynomial mapping has the form $R_{Nn}(x) =P_n(\pi_N(x))$, where $\pi_N(x)$ is a polynomial of degree $N$. Again it is assumed that both $P_n(x)$ and $R_n(x)$ are nondegenerate orthogonal polynomials. The theory of such mappings was considered in \cite{GVA}.

Consider now the following concrete example connected with Jacobi polynomials.
This example will allow to establish the weight function of the big $-1$-Jacobi polynomials.

Let
$$
P_n^{(\xi,\eta)}(x) = G_n \: {_2}F_1 \left( {-n, n+\xi+\eta+1 \atop \xi+1}; x \right)
$$
be Jacobi polynomials with orthogonality relation
$$
\int_0^1 x^{\xi} (1-x)^{\eta} P_n^{(\xi,\eta)}(x) P_m^{(\xi,\eta)}(x) dx = h_nm \: \delta_{nm}
$$
on the interval $[0,1]$.

The normalization coefficient
$$
G_n = (-1)^n \: \frac{(\xi+1)_n}{(n+\xi+\eta+1)_n}
$$
ensures that $P_n(x)$ is monic $P_n^{(\xi,\eta)}(x)=x^n + O(x^{n-1})$.

Perform first an affine transformation of the argument and consider the new monic orthogonal polynomials
$$
P_n(x)= (c^2-1)^n \: P_n^{(\xi,\eta)}\left(\frac{1-x}{1-c^2}\right),
$$
where $c$ is a real parameter with the restriction $0<c<1$.

In terms of hypergeometric functions
\be
P_n(x) = (1-c^2)^n \: \frac{(\xi+1)_n}{(n+\xi+\eta+1)_n} \: {_2}F_1 \left( {-n, n+\xi+\eta+1 \atop \xi+1}; \frac{1-x}{1-c^2} \right) \lab{hyp_P} \ee
Clearly, these polynomials are orthogonal on the interval $[c^2,1]$
$$
\int_{c^2}^1 (1-x)^{\xi} (x-c^2)^{\eta} P_n(x) P_m(x) dx = 0, \quad n \ne m
$$
Introduce also the companion polynomials $Q_n(x)$ through the Christoffel transform
\be
Q_n(x) = \frac{P_{n+1}(x) - A_n P_n(x)}{x-1}, \quad A_n=\frac{P_{n+1}(1)}{P_n(1)} \lab{QP_Jac} \ee
It is easily seen that the polynomials $Q_n(x)$ are again expressible in terms of Jacobi polynomials with $\xi \to \xi+1$:
$$
Q_n(x) = (c^2-1)^n \:  \: P_n^{(\xi+1,\eta)}\left(\frac{1-x}{1-c^2}\right),
$$
or, in terms of hypergeometric functions
\be
Q_n(x) =(1-c^2)^n \: \frac{(\xi+2)_n}{(n+\xi+\eta+2)_n} \: {_2}F_1 \left( {-n, n+\xi+\eta+2 \atop \xi+2}; \frac{1-x}{1-c^2}\right). \lab{hyp_Q} \ee
The polynomials $P_n(x)$ and $Q_n(x)$ are connected by the relations \re{Q_P_CT} and \re{PQ_GT} with $\nu=1$.
The coefficients $A_n$ and  $B_n$ can be found from the following observation. Putting $x=1$, we find from \re{hyp_P} and \re{hyp_Q}
$$
P_n(1) = (1-c^2)^n \: \frac{(\xi+1)_n}{(n+\xi+\eta+1)_n}, \quad Q_n(1) = (1-c^2)^n \: \frac{(\xi+2)_n}{(n+\xi+\eta+2)_n}.
$$
From these formulas we immediately get
\be
A_n= \frac{P_{n+1}(1)}{P_n(1)} = (1-c^2) \: \frac{(\xi+n+1)(\xi+\eta+n+1)}{(2n+\xi+\eta+1)(2n+\xi+\eta+2)} \lab{A_n_Jac} \ee
and
\be
B_n= \frac{Q_n(1)-P_n(1)}{Q_{n-1}(1)} = (1-c^2) \: \frac{n(\eta+n)}{(2n+\xi+\eta)(2n+\xi+\eta+1)}. \lab{B_n_Jac} \ee
Note that $A_n>0, \: B_n>0$ for $n=1,2,\dots$ due to the restriction $0<c<1$.

Consider now the new monic orthogonal polynomials $R_n(x)$ defined by the relations
\be
R_{2n}(x) = P_n(x^2), \quad R_{2n+1}(x) = (x-1) Q_n(x^2) \lab{R_PQ_Jac} \ee
According to the general theory of polynomial mappings \cite{GVA}, \cite{MP}, it is not difficult to show that the polynomials $R_n(x)$
are orthogonal on a domain formed by the union of {\it two} intervals $[-1,-c], \: [c,1]$ of the real axis:
\be
\int_{-1}^{-c} R_n(x) R_m(x) W(x) dx + \int_{c}^{1} R_n(x) R_m(x) W(x) dx = 0, \quad n \ne m \lab{ort_R_Jac} \ee
where the (non-normalized) weight function is:
\be
W(x) = \theta(x)(1+x) (1-x^2)^{\xi} (x^2-c^2)^{\eta} \lab{W_R_Jac} \ee
and $\theta(x)=x/|x|$ is the sign function. Note that the weight function $W(x)$ is not positive on the interval $[-1,-c]$.

In terms of Gauss' hypergeometric functions we have the expressions
\be
R_{2n}(x) = P_n(x^2) = (1-c^2)^n \: \frac{(\xi+1)_n}{(n+\xi+\eta+1)_n} \: {_2}F_1 \left( {-n, n+\xi+\eta+1 \atop \xi+1}; \frac{1-x^2}{1-c^2} \right) \lab{R_2n_hyp} \ee
and
\be
R_{2n+1}(x) = (x-1)Q_n(x^2) = (1-c^2)^n \: \frac{(\xi+2)_n}{(n+\xi+\eta+2)_n} \: (x-1){_2}F_1 \left( {-n, n+\xi+\eta+2 \atop \xi+2}; \frac{1-x^2}{1-c^2} \right) \lab{R_2n+1_hyp} \ee
The polynomials $R_n(x)$ satisfy the 3-term recurrence relation
\be
R_{n+1}(x) + (-1)^n R_n(x) +  v_n R_{n-1}(x)=xR_n(x), \lab{R_Jac_rec} \ee
where
\be
v_{2n}=-B_n= (c^2-1) \: \frac{n(\eta+n)}{(2n+\xi+\eta)(2n+\xi+\eta+1)}, \quad v_{2n+1}=-A_n=(c^2-1) \: \frac{(\xi+n+1)(\xi+\eta+n+1)}{(2n+\xi+\eta+1)(2n+\xi+\eta+2)}
\lab{v_n_expr} \ee
Note that all the coefficients $v_n$ are negative, $v_n<0$, which corresponds to the non-positivity of the weight function $W(x)$.

\section{The weight function and the orthogonality of the big -1 Jacobi polynomials}
\setcounter{equation}{0}
In order to determine the weight function and the orthogonality region for the big -1 Jacobi polynomials,
we notice that formulas \re{even_hyp} and \re{odd_hyp} can be presented in the following equivalent form:
\be
P_n^{(-1)}(x) = R_n(x) - G_n R_{n-1}(x), \lab{Ger_PR} \ee
where $R_n(x)$ are the polynomials defined by \re{R_2n_hyp}, \re{R_2n+1_hyp} that satisfy the recurrence relation \re{R_Jac_rec}.
In these formulas we should put $\xi=(\alpha-1)/2, \; \eta=(\beta+1)/2$.   The coefficients $G_n$ have the expression
\be
G_n = \left\{ {  \frac{(1-c)n}{2n+\alpha+\beta},  \quad n \quad \mbox{even}   \atop  -\frac{(1+c)(n+\alpha)}{2n+\alpha+\beta},  \quad n \quad \mbox{odd}}  \right . \lab{G_n_expr} \ee
It is well known that if two families of orthogonal polynomials are related by a formula such as \re{Ger_PR},
then necessarily the polynomials $P_n^{(-1)}(x)$ are obtained from the polynomials $R_n(x)$ by the Geronimus transform \cite{ZhS}.
This is equivalent to the statement that the weight function $w^{(-1)}(x)$ of the polynomials $P_n^{(-1)}(x)$
is obtained from the weight function $W(x)$ of the polynomials  $R_n(x)$ as follows:
\be
w^{(-1)}(x) = \frac{W(x)}{x-\mu} + M \delta(x-\mu), \lab{ww_ger} \ee
with two additional parameters $\mu$ and $M$. Formula \re{ww_ger} means that apart from the division of the weight function $W(x)$
by the linear factor $x-\mu$ there is an additional concentrated mass $M$ that is inserted at the point $x= \mu$.

The parameter $\mu$ can be found from the recurrence relation for the coefficients $G_n$ \cite{ZhS}
\be
G_{n+1} +(-1)^n + \frac{v_n}{G_n} = \mu \lab{rec_G} \ee
with the recurrence coefficients $v_n$ given by \re{v_n_expr}.

Substituting \re{G_n_expr} into \re{rec_G} we obtain $\mu=-c$.

Thus the orthogonality relation for polynomials $P_n^{(-1)}(x)$ takes the form
\be
\int_{\Gamma} P_n^{(-1)}(x)P_m^{(-1)}(x) W(x)(x+c)^{-1} dx + M P_n^{(-1)}(-c)P_m^{(-1)}(-c) =0, \quad n \ne m \lab{ort_P_-1_M} \ee
where the contour $\Gamma$ is the union of the two intervals $[-1,-c]$ and $[c,1]$ of the real axis.

In order to find the value $M$ of the concentrated mass it is sufficient to consider a special case
of the orthogonality relation  \re{ort_P_-1_M} for $n=1,m=0$
\be
\int_{\Gamma} P_1^{(-1)}(x) W(x) (x+c)^{-1} dx +M P_1(-c) =0 \lab{ort_1_pol} \ee
Now, $P_1^{(-1)}(x)$ is given by
$$
P_1^{(-1)}(x)=x+\zeta
$$
where $\zeta = \frac{c(\alpha+1)-\beta-1}{2+\alpha+\beta}$. Substituting this expression into \re{ort_1_pol}
and calculating the integral (through an elementary reduction to the Euler beta-integral) we find that $M=0$.

Thus, the orthogonality relation for polynomials $P_n^{(-1)}(x)$ reads
\be
\int_{\Gamma} P_n^{(-1)}(x)P_m^{(-1)}(x) w^{(-1)}(x) dx = 0, \quad n \ne m, \lab{ort_P-1} \ee
where the weight function $w^{(-1)}(x)$ can be presented in the form
\be
w^{(-1)}(x) = \theta(x) (x+1)(x+c)^{-1} (1-x^2)^{(\alpha-1)/2}(x^2 -c^2)^{(\beta+1)/2} \lab{w-1} \ee
or, equivalently
\be
w^{(-1)}(x) = \theta(x) (x+1)(x-c) (1-x^2)^{(\alpha-1)/2}(x^2 -c^2)^{(\beta-1)/2}. \lab{w-12} \ee
Note that under the restrictions $\alpha>-1,\beta>-1$, the weight function is positive on the two intervals of $\Gamma$ and all the moments
$$
m_n = \int_{\Gamma} w^{(-1)}(x) x^n dx
$$
are finite for $n=0,1,2,\dots$.

When $c \to 0$, the big q-Jacobi polynomials reduce to the little q-Jacobi polynomials \cite{KS}. In the limit case $q \to -1$ we see that when
$c \to 0$, the set of two intervals coalesces to the single interval $[-1,1]$ and the weight function becomes
$$
w(x)\left|_{c=0} \right . = (1+x) |x|^{\beta}(1-x^2)^{(\alpha-1)/2}
$$
which corresponds to the weight function of the little $-1$ Jacobi polynomials \cite{VZ_little}.

So far, we considered the case $0<c<0$. The case $c>1$ can be treated in an analogous way. This leads to the following orthogonality relation
\be
\int_{-c}^{-1}P_n^{(-1)}(x) P_m^{(-1)}(x) w^{(-1)}(x) dx + \int_{1}^{c}P_n^{(-1)}(x) P_m^{(-1)}(x) w^{(-1)}(x) dx =0,\quad n \ne m, \lab{c>1_ort} \ee
where the weight function is almost the same as in \re{w-12} with obvious modifications:
\be
w^{(-1)}(x) = \theta(x) (x+1)(c-x) (x^2-1)^{(\alpha-1)/2}(c^2-x^2)^{(\beta-1)/2} \lab{wc>1} \ee
and where again we have the restrictions $\alpha>-1, \: \beta>-1$.

The case $c=1$ is degenerate: the recurrence coefficients $u_{n}$ for even $n$ become zero, $u_{2n}=0$,
which means that orthogonal polynomials $P_n^{(1)}(x)$ are no more positive definite.
The two intervals of orthogonality shrink  into two points $x=\pm 1$.

\section{The big -1 Jacobi polynomials as limits of the Bannai-Ito polynomials}
\setcounter{equation}{0} Bannai and Ito proved a theorem which
characterizes all dual polynomial schemes \cite{BI}. This theorem
is a generalization of the Leonard theorem \cite{Leonard}, where
only the finite-dimensional case was considered. Apart from the
already known explicit orthogonal polynomials from the Askey
scheme, Bannai and Ito found one more example corresponding to the
limit $q \to -1$ of the q-Racah polynomials. The Bannai-Ito
polynomials $W_n(x)$ are orthogonal on a finite number of points
$x_i, \; i=0,1,2,\dots,N$ given by \be x_i = \left\{ { \theta_0
+2hi, \quad \mbox{} \quad i \quad \mbox{even}  \atop  \theta_0 -
2h(i+1-s), \quad \mbox{} \quad i \quad \mbox{odd}}  \right . ,
\lab{x_BI} \ee where $\theta_0, h, s$ are arbitrary parameters.

The Bannai-Ito polynomials satisfy the recurrence relation
\cite{BI}, \cite{Ter2} \be A_n W_{n+1}(x) +(\theta_0 -A_n -
C_n)W_n(x) +C_n W_{n-1}(x) =x W_n(x), \lab{rec_BI} \ee where \be
A_n=\left\{ { \frac{2h(n+1+r_1)(n+1+r_2)}{2n+2-s^*}, \quad \mbox{}
\quad n \quad \mbox{even} \atop
\frac{2h(n+1-s^*)(n+1-r_3)}{2n+2-s^*}  \quad \mbox{} \quad n \quad
\mbox{odd}} \right . \lab{A_BI} \ee and \be C_n=\left\{ {
-\frac{2h n(n-s^*+r_3)}{2n-s^*}, \quad \mbox{} \quad n \quad
\mbox{even} \atop -\frac{2h(n-r_1-s^*)(n-r_2-s^*)}{2n-s^*} \quad
\mbox{} \quad n \quad \mbox{odd}} \right . ,\lab{C_BI} \ee where
$s^*, r_1,r_2$ are parameters such that \be s+s^*=r_3-r_1-r_2
\lab{ss_rr} \ee and $r_2=-N-1$ if $N$ is even and $r_3=N+1$ if $N$
is odd.

The polynomials $W_n(x)$ are not monic. Instead, they satisfy the
conditions \cite{BI} $W_n(\theta_0) =1, \: n=0,1,2,\dots,N$. As
usual, it is assumed that $W_{-1}(x)=0$. We can introduce the
monic Bannai-Ito polynomials $\hat W_n(x)=x^n + O(x^{n-1})$ by the
formula
$$\hat W_n(x) = A_0 A_1 \dots A_{n-1} W_n(x).$$
The polynomials $\hat W_n(x)$ are then seen to satisfy the
recurrence relation \be \hat W_{n+1}(x) + (\theta_0 -A_n-C_n) \hat
W_{n}(x) + A_{n-1} C_n \hat W_{n-1}(x) = x \hat W_{n}(x).
\lab{rec_mW} \ee Note that the parameter $\theta_0$ is not
essential: it only specifies the initial point $x_0$. We can thus
put $\theta_0=1$ without loss of generality. The recurrence
relation \re{rec_mW} then determines all polynomials $\hat W_n(x)$
uniquely.

Consider now a limit $N \to \infty$. Given the relations that
$r_2$ and $r_3$ have with $N$, it is clear that under such a
limit, we must necessarily have: $r_3 \to \infty$ and $r_2 \to
-\infty$. We assume that the parameters $r_1, s^*$ are finite. It
is now convenient to put
$$
r_1 = \alpha, \quad s^*=-\alpha-\beta
$$
with some fixed parameters $\alpha, \beta$.

From \re{ss_rr} it then follows that $s \to \infty$. The parameter
$h$ moreover, is free: it allows for a scaling of the argument
$x$. We can thus choose $h \to 0$ in such a way that \be 2h r_2
\to c+1, \quad 2h r_3 \to c-1, \lab{r_23} \ee where $c$ is a
parameter such that $0<c<1$.

Then, in the limit $N \to \infty$, the recurrence coefficients
$A_n$ and $C_n$ coincide with the coefficients $A_n,C_n$ given by
formulas \re{A_lim} and \re{C_lim}, and we have already observed
that these coefficients $A_n$ and $C_n$ determine uniquely the big
-1 Jacobi polynomials.

Consider the behavior of the argument $x_i$ in the limit $N \to
\infty$. As $h \to 0$, let $i \to \infty$ in a way such that $2hi
\to -y$, where $y$ is a continuous variable which varies from 0 to
$c$. Hence from \re{x_BI}, in the limit $N \to \infty$,  the
points $x_i$, for even $i$, densely cover the interval $[c,1]$.
Similarly, from \re{ss_rr}, \re{r_23}
$$
\lim_{N \to \infty} 2hs = \lim_{N \to \infty} \{2hr_3 - 2h r_2\}
=-2
$$
and hence from \re{x_BI}, in the limit $N\to \infty$, the points
$x_i$, for odd $i$, densely cover the interval $[-1,-c]$. This
corresponds to the fact that the big -1 Jacobi polynomials are
orthogonal on the union of the two intervals $[-1,-c]$ and
$[c,1]$.

We thus have shown that the big -1 Jacobi polynomials can be
obtained from the Bannai-Ito polynomials $W_n(x)$ through a limit
process where $N \to \infty$. This complements the fact
demonstrated in this paper
 that the big -1 Jacobi polynomials arise as a $q \to -1$ limit
from the big q-Jacobi polynomials.

When we remind ourselves that the big q-Jacobi polynomial result
from a $N \to \infty$ limit of the q-Racah polynomials \cite{KS},
we see that we have the following commutative diagram for limit
relations between these polynomials:

\vspace{8mm}

\begin{displaymath}
    \xymatrix{
        \boxed{\mbox{q-Racah}} \ar[r]^{q \to -1} \ar[d]_{N \to \infty} & \boxed{\mbox{Bannai-Ito}} \ar[d]^{N \to \infty} \\
        \boxed{\mbox{Big q-Jacobi}} \ar[r]_{q \to -1}       & \boxed{\mbox{Big  -1 Jacobi}} }
\end{displaymath}

\vspace{8mm}

\section{Anticommutator algebra describing big -1 Jacobi polynomials}
\setcounter{equation}{0}
The Askey-Wilson polynomials are described by the AW(3)-algebra \cite{Zhe}, \cite{Ter}.
Among the different equivalent forms of this algebra, we choose the following one:
\be
XY-qYX=\mu_3 Z+\omega_3, \quad YZ-qZY=\mu_1 X+\omega_1, \quad ZX-qXZ=\mu_2 Y+\omega_2, \lab{AW3} \ee
which possesses an obvious symmetry with respect to all 3 operators (see, e.g. \cite{IT}).

Here $q$ is a fixed parameter corresponding to the "base" parameter in the q-hypergeometric functions defining the Askey-Wilson polynomials \cite{KS}.
The pairs of operators $(X,Y)$, $(Y,Z)$ and $(Z,X)$ play the role of "Leonard pairs" (see \cite{Ter}, \cite{IT}).

The Casimir operator
\be
Q= (q^2-1)XYZ + \mu_1 X^2 + \mu_2
q^2 Y^2 + \mu_3 Z^2 +(q+1)(\omega_1 X + \omega_2 q Y + \omega_3 Z) \lab{Cas_Q} \ee
commutes with all operators $X,Y,Z$.

The constants $\omega_i, \: i=1,2,2$ (together with the value of the Casimir operator $Q$)
define representations of the $AW(3)$ algebra (see \cite{Zhe} for details).

Consider now the case of  the big $-1$ Jacobi polynomials and choose the following operators
\be
X=L_0+\alpha+\beta+1, \quad Y=x, \quad Z= -\frac{2}{x}\left(c+(x-1)(x+c)R  \right),
\lab{XYZ_real} \ee
where $L_0$ is the operator given by \re{L_0}.

It is then easy to verify that these operators satisfy the linear anticommutation relations
\be
XY+YX = Z+\omega_3, \quad YZ+ZY=\omega_1, \quad ZX+XZ=4Y+\omega_2, \lab{anti_alg} \ee
where
$$
\omega_1=-4c, \; \omega_2 = 4(\alpha-\beta c), \; \omega_3 =2(\beta-\alpha c).
$$
The Casimir operator of the algebra defined by \re{anti_alg} is
$$
Q=Z^2+4Y^2
$$
In the realization \re{XYZ_real} the Casimir operator takes the constant value $Q=4(c^2+1)$.

In this realization the operator $X$ (up to an additive constant) is the operator of which the polynomials $P_n^{(-1)}(x)$ are the eigenfunctions.
The operator $Y$ here corresponds to multiplication by $x$.

The ``dual''realization of the algebra \re{anti_alg} is obtained
if one takes an infinite discrete basis $e_n, \; n=0,1,2,\dots$ on
which the operators $X,Y$ act as \be X e_n=(\lambda_n +
\alpha+\beta+1)e_n, \quad Y e_n  = u_{n+1}^{(-1)}e_{n+1} +
b_n^{(-1)} e_n + e_{n-1}, \lab{XY_dual} \ee where $\lambda_n$ is
the eigenvalue \re{lam-1} and where the recurrence coefficients
$u_n^{(-1)}, b_n^{(-1)}$  are given by \re{u_lim}, \re{b_lim}.
Thus in this representation the operator $Y$ is a Jacobi (i.e.
tri-diagonal) matrix and the eigenvalue equation
$$
Y \vec P = x \vec P
$$
is equivalent to the recurrence relation \re{rec-1} for the big $-1$ Jacobi polynomials.
Indeed, we can present the vector $\vec P$ in terms of its expansion coefficients over the basis $e_n$:
$$
\vec P = \sum_{n=0}^{\infty} C_n e_n
$$
Without loss of generality we can choose $C_0=1$.
The coefficients $C_n$ in this expansion are then found to satisfy the recurrence relation \re{rec-1}
and  it is seen moreover that these $C_n$ are monic polynomials in $x$ of degree $n$. Hence $C_n=P_n^{(-1)}(x)$.

\section{Two-diagonal basis for the operator $L_0$ and a generalization of Gauss' hypergeometric functions}
\setcounter{equation}{0}
We already showed that the Dunkl-type operator $L_0$ is tri-diagonal in the ordinary monomial basis $x^n$  (see formula \re{L-3-diag_x}).
There exists, however, another polynomial basis in which the operator $L_0$ is two-diagonal. This basis can be constructed as follows
\be \phi_0=1, \: \phi_1(x)=x-1, \phi_2(x) =(x^2-1), \dots, \phi_{2n}(x)=(x^2-1)^n, \: \phi_{2n+1}(x) =(x-1)(x^2-1)^n \lab{phi_basis} \ee
It is easily verified that
\be
L_0 \phi_n(x) = \lambda_n \phi_n(x) + \eta_{n} \phi_{n-1}(x), \lab{L_0_phi} \ee
where $\lambda_n$ is the eigenvalue given by by \re{lam-1}, and
$$
\eta_n = \left\{ {2n(c-1), \quad \mbox{if} \quad n \quad \mbox{even}  \atop  -2(c+1)(\alpha+n) \quad \mbox{if} \quad n \quad \mbox{odd}}  \right .
$$
Consider now the eigenvalue equation
\be
L_0 P_n(x) = \lambda_n P_n(x)
\lab{L0PP} \ee
and expand the polynomials $P_n(x)$ over the basis $\phi_n(x)$:
$$
P_n(x) = \sum_{s=0}^n A_{ns} \phi_s(x).
$$
For the expansion coefficients $A_{ns}$ we have from \re{L0PP}:
\be
A_{n,s+1} = \frac{A_{ns}(\lambda_n-\lambda_s)}{\eta_{s+1}} \lab{rec_A_s} \ee
From \re{rec_A_s}, the coefficients $A_{ns}$ can be found explicitly in terms of $A_{n0}$:
\be A_{ns}= A_{n0} \: \frac{(\lambda_n-\lambda_0)(\lambda_n-\lambda_1) \dots (\lambda_n-\lambda_{s-1})}{\eta_1 \eta_2 \dots \eta_s} \lab{A_ns_A_0} \ee
or in terms of the coefficient $A_{nn}$:
\be
A_{ns}= A_{nn} \: \frac{\eta_n \eta_{n-1} \dots \eta_{s+1}}{(\lambda_n-\lambda_{n-1})(\lambda_n - \lambda_{n-2}) \dots (\lambda_n - \lambda_s) } \lab{A_ns_A_n} \ee
We thus have the following explicit formula for the polynomials $P_n(x)$
\be P_n(x) = A_{n0} \: \sum_{s=0}^n \frac{(\lambda_n-\lambda_0)(\lambda_n-\lambda_1) \dots (\lambda_n-\lambda_{s-1})}{\eta_1 \eta_2 \dots \eta_s} \phi_s(x) \lab{P_n_hyp_phi} \ee
Expression \re{P_n_hyp_phi} resembles Gauss' hypergeometric function and can be considered as a nontrivial generalization of it.
Indeed, products in the numerator and denominator of \re{P_n_hyp_phi} can easily be presented in terms of ordinary Pochhammer symbols
and we thus recover the explicit formulas \re{even_hyp} and \re{odd_hyp}. Note nevertheless, that the form \re{P_n_hyp_phi} looks much simpler.

Moreover, note also that in the basis $\phi_n(x)$ the operators $X$ and $Y$ of the algebra defined by \re{XYZ_real} and \re{anti_alg},
become lower and upper triangular:
$$
X \phi_n(x) = (L_0+\alpha+\beta+1)\phi_n(x) = (\lambda_n+\alpha+\beta+1) \phi_n(x) + \eta_n \phi_{n-1}(x)
$$
and
$$
Y \phi_n(x) = x \phi_n(x) = \phi_{n+1}(x)+(-1)^n \phi_{n}(x)
$$
Note that such bases in which the operators $X,Y$ are
two-diagonal, were central objects in Terwilliger's approach to
Leonard pairs \cite{Ter}, \cite{Ter2}. Formulas similar
to \re{P_n_hyp_phi} also appear in the theory of Leonard pairs.

\bigskip\bigskip
{\Large\bf Acknowledgments}
\bigskip

\noindent The authors are indebted for R.Askey, C.Dunkl, T.
Koornwinder, W.Miller, P.Terwilliger and P.Winternitz for stimulating communications, AZ
thanks CRM (U de Montr\'eal) for its hospitality.

\newpage

\bb{99}

\bi{AA} G.E.Andrews and R.Askey, {\it Classical orthogonal
polynomials}, Polyn\^omes  Orthogonaux et Applications, Lecture
Notes in Mathematics, 1985, V. {\bf 1171}, 36--62.

%\bi{biask} R. Askey, Comments in: "Gabor Szeg\H{o}:Collected Papers",
%Birkh\"auser, Basel, 1982, v.1, 303-305.

\bi{AI2} R. A. Askey and M. E. H. Ismail, {\it A generalization of
ultraspherical polynomials}, Studies in Pure Mathematics (P.
Erd\H{o}s, ed.), Birkhauser, Basel, 1983, pp. 55-78.

\bi{AI} R. A. Askey and M. E. H. Ismail, {\it Recurrence
relations, continued fractions and orthogonal polynomials}, Mem.
AMS, {\bf 49} (1984), No. 300,  1--108.

%\bi{Atia} M. J. Atia, {\it An example of non-symmetric
%semi-classical form of class s=1. Generalization of a case of
%Jacobi sequence}. Int. J. Math. and Math. Sci. {\bf 24} (10)
%(2000), 673--689.

%\bi{Atia2} M.J.Atia, {\it Some generalized Jacobi polynomials},
%OPSFA 2009, http://wis.kuleuven.be/OPSFA/bookletOPSFA.pdf

%\bi{Bax} G. Baxter, {\it Polynomials defined by a difference
%system}, J.Math.Anal.Appl. {\bf 2} (1961), 223-263.

\bi{BI} E. Bannai and T. Ito, {\it Algebraic Combinatorics I:
Association Schemes}. 1984. Benjamin \& Cummings, Mento Park, CA.

\bi{Cheikh} Y. Ben Cheikh and M.Gaied, {\it Characterization of the Dunkl-classical symmetric orthogonal polynomials}, Appl. Math. and Comput.
{\bf 187}, (2007) 105--114.

%\bi{Chi2} T.Chihara, {\it Orthogonal polynomials with Brenke type
%generating functions}, Duke Math. J. {\bf 35} (1968), 505--517.

\bi{Chi} T. Chihara, {\it An Introduction to Orthogonal
Polynomials}, Gordon and Breach, NY, 1978.

\bi{Chi_Chi} L.M.Chihara,T.S.Chihara, {\it A class of nonsymmetric orthogonal polynomials}. J. Math. Anal. Appl. {\bf 126} (1987), 275--291.

\bi{Dunkl} C.F.Dunkl, {\it Integral kernels with reflection group invariance}. Canadian Journal of Mathematics, {\bf 43} (1991)
1213--1227.

%\bi{Ger} Ya.L. Geronimus,\quad {\it Polynomials Orthogonal on a Circle
%and their Applications}, \\ Am.Math.Transl.,Ser.1 {\bf 3}(1962), 1-78.

\bi{Ger1} Ya.L.Geronimus, {\it On polynomials orthogonal with respect to to
the given numerical sequence and on Hahn's theorem}, Izv.Akad.Nauk, {\bf 4}
(1940), 215--228 (in Russian).

\bi{GVA} J.Geronimo and W. Van Assche, {\it Orthogonal polynomials on several intervals via a polynomial mapping},
Trans. Amer. Math. Soc., {\bf 308}(2) (1988), 559--581.

%\bi{Ger2} J. Geronimus, {\it The orthogonality of some systems of
%polynomials}, Duke Math. J. {\bf 14} (1947), 503--510.

\bi{Hahn} W. Hahn, \"Uber Orthogonalpolynome die
q-Differenzengleichungen gen\"ugen, Math. Nath. {\bf 2} (1949),
4--34.

%\bi{Hahn} W.Hahn, {\it \"Uber die Jacobischen Polynome und Zwei
%verwandte Polynomklassen}, Math.Z. {\bf 39} (1935), 634-638.

%\bi{HN} E.Hendriksen and O. Nj{\aa}stad, {\it Biorthogonal Laurent polynomials
%with biorthogonal derivatives}, Rocky Mount. J. Math. {\bf 21} (1991),
%391-317.

%\bi{HR} E.Hendriksen and H. van Rossum, {\it Orthogonal Laurent polynomials},
%Indag. Math. (Ser. A) {\bf 48} (1986), 17-36.

%\bi{IsMas} M.E.H. Ismail and D. Masson, {\it Generalized orthogonality and
%continued fractions}, \\J.Approx.Theory {\bf 83} (1996), 1-40.

%\bi{JT} W.B.Jones and W.J.Thron, {\it Survey of continued fraction methods of
%solving moment problems} in: analytic Theory of Continued Fractions, LNM 932,
%Springer, Berlin, Heidelberg, New York (1981).

%\bi{KLY} K.H. Kwon, L.L. Littlejohn and B.H. Yoo, {\it New characterizations of classical orthogonal polynomials}, Indag. Mathem., N.S., {\bf 7} (1996), %199--213.

\bi{IT} T.Ito, P.Terwilliger, {\it Double Affine Hecke Algebras of Rank 1
and the $Z_3$-Symmetric Askey-Wilson Relations}, SIGMA {\bf 6} (2010), 065, 9 pages.

\bi{KS} R.Koekoek, R.Swarttouw, {\it The Askey-scheme of hypergeometric orthogonal polynomials and its q-analogue}, Report no. 98-17, Delft University of Technology, 1998.

\bi{KLS} R. Koekoek,P. Lesky, R. Swarttouw, {\it Hypergeometric Orthogonal Polynomials and Their Q-analogues}, Springer-Verlag, 2010.

%\bi{Kwon} K.H. Kwon, L.L. Littlejohn, B.H. Yoo, {\it
%Characterization of orthogonal polynomials satisfying differential
%equations}, SIAM J. Math. Anal. {\bf 25} (1994) 976-�990.

\bi{Leonard} D.Leonard, {\it Orthogonal Polynomials, Duality and
Association Schemes}, SIAM J. Math. Anal. {\bf 13} (1982)
656--663.

%\bi{Mag} A.P. Magnus, {\it Painleve\'e-type differential equations
%for the recurrence coefficients of semi-classical orthogonal
%polynomials}, J. Comp. Appl. Math. {\bf 57} (1995), 215-237.

\bi{MP} F.Marcell\'an and J.Petronilho, {\it Eigenproblems for Tridiagonal 2-Toeplitz Matrices
and Quadratic Polynomial Mappings},  Lin. Alg. Appl. {\bf 260} (1997)   169--208.

%\bi{Mar} P.Maroni, {\it Variations around classical orthogonal polynomials.
%Connected problems}, J.Comp.Appl.Math. {\bf 48} (1993), 133-155.

%\bi{Mass} D.Masson, {\it Difference equations, continued fractions, Jacobi
%matrices and orthogonal polynomials}, in Nonlinear Numerical methods and
%Rational Approximations, (A.Cuyt, ed.), 239-257. Reidel Publ.Co. 1988.

%\bi{NSU} A.F. Nikiforov, S.K. Suslov, and V.B. Uvarov, {\em
%Classical Orthogonal Polynomials of a Discrete Variable},
%Springer, Berlin, 1991.

%\bi{Pas} P.I. Pastro, {\it Orthogonal polynomials and some q-beta
%integrals of Ramanujan},\\  J.Math.Anal.Appl. {\bf 112} (1985),
%517--540.

%\bi{SVZ} V.Spiridonov, L.Vinet and A.Zhedanov, {\it Spectral
%transformations, self-similar reductions and orthogonal polynomials}, J.Phys.
%A:  Math.  and Gen.  {\bf 30} (1997), 7621--7637.

\bi{Sz} G. Szeg\H{o}, Orthogonal Polynomials, fourth edition,  AMS, 1975.

\bi{Ter} P. Terwilliger, {\it Two linear transformations each tridiagonal with respect to an eigenbasis of the other}.
Linear Algebra Appl. {\bf 330} (2001) 149--203.

\bi{Ter2} P.Terwilliger, {\it Two linear transformations each
tridiagonal with respect to an eigenbasis of the other; an
algebraic approach to the Askey scheme of orthogonal polynomials},
arXiv:math/0408390v3.

%\bi{VYZ} L.Vinet, O.Yermolayeva and A.Zhedanov, {\it A method to
%study the Krall and q-Krall polynomials}, J.Comp.Appl.Math. {\bf
%133} (2001) 647--656.

%\bi{VZ_Bochner} L.Vinet and A.Zhedanov, {\it Generalized Bochner
%theorem: characterization of the Askey-Wilson polynomials},
%J.Comp.Appl.Math., {\bf 211} (2008) 45 -- 56.

\bi{VZ_little} L.Vinet and A.Zhedanov, {\it A ``missing`` family
of classical orthogonal polynomials}, arXiv:1011.1669v2.

%\bi{WW} E.T. Whittacker, G.N. Watson, {\em A Course of Modern
%Analysis}, Cambridge, 1927.

\bi{Zhe} A. S. Zhedanov. {\it Hidden symmetry of Askey-Wilson polynomials}, Teoret. Mat. Fiz.
{\bf 89} (1991) 190--204. (English transl.: Theoret. and Math. Phys. {\bf 89} (1991), 1146--1157).

\bibitem{ZhS} A.S. Zhedanov, {\it Rational spectral transformations
and orthogonal polynomials}, J. Comput. Appl. Math. {\bf 85}, no. 1
(1997), 67--86.

\eb

\end{document}